\documentclass[10pt,twocolumn]{article}

\usepackage{fullpage}
\usepackage{amssymb}
\usepackage{xcolor}
\usepackage{mathtools}

\newtheorem{theorem}{Theorem}

\newtheorem{assumption} {Assumption}

\newtheorem{lemma} {Lemma}

\newtheorem{example}{Example}
\newtheorem{definition} {Definition}
\def\ds{\displaystyle}
\def\s{\bar t}
\def\Rset{{\rm I}\!{\rm R}}
\def\qed{\hfill \ensuremath{\Box}}
\def\ds{\displaystyle}
\def\s{\sigma}
\def\S{\Sigma}
\def\R{\mathbb{R}}

\def\C{\mathcal{C}}
\def\U{\mathcal{U}}

\def\Rset{{\rm I}\!{\rm R}}
\def\qed{\hfill \ensuremath{\Box}}

\title{\LARGE \bf ISS characterization of retarded switching systems with relaxed Lyapunov--Krasovskii functionals
\footnote{This work is supported by a public grant overseen by the Research and Valorization Service (SRV) of ENSEA.}}

\author{Ihab Haidar\thanks{Quartz EA 7393, ENSEA, Cergy-Pontoise, France, {\tt ihab.haidar@ensea.fr}} and 
Pierdomenico Pepe\thanks{Department of Information Engineering, Computer Science, and Mathematics, University of L'Aquila, 67100 L'Aquila, Italy {\tt pierdomenico.pepe@univaq.it}.}}

\begin{document}

\maketitle
\thispagestyle{empty}
\pagestyle{empty}

\begin{abstract}
This paper gives further insights about the Lyapunov--Krasovskii characterization of input-to-state stability (ISS) for switching retarded systems on the basis of the results in {\it [I. Haidar and P. Pepe. Lyapunov–krasovskii characterization of the input-to-state stability for switching retarded systems. SIAM
Journal on Control and Optimization, 59(4):2997–3016, 2021]}. We give new characterizations of the ISS property through the existence of a relaxed common Lyapunov-Krasovskii functional. More precisely, we show that the existence of a continuous Lyapunov-Krasovskii functional whose upper right-hand Dini derivative satisfies a dissipation inequality almost everywhere is necessary and sufficient for the ISS of switching retarded systems with measurable inputs and measurable switching signals. Different characterization results, using different derivative notions, are also given. 
\end{abstract}

\textbf{Keywords:}
Input-to-state stability; Converse theorems; Lyapunov--Krasovskii functionals; retarded functional differential equations; switching systems.

\section{Introduction}
The problem of stability of switching systems has attracted much attention in the literature of control theory (see, e.g., \cite{AL2001, Boscain2002, Liberzon, Liberzon-Morse, Mazenc2018, SSGE, WANG201678, YanOzbay2008} and the references therein). The existence of a common Lyapunov function, i.e., a function which decreases uniformly along the trajectories of individual subsystems, is a sufficient condition for various stability notions like uniform asymptotic, exponential, and input-to-state stability. The existence of a common Lyapunov function is also necessary for switching systems which are uniformly stable. Converse Lyapunov theorems characterizing the stability of a switching system by the existence of a common Lyapunov function have been then developed for various switching dynamics (see, e.g.,~\cite{Dayawansa,MANCILLAAGUILAR200067,Wirth2005}) for finite-dimensional systems,~\cite{haidar2021lyapunov, Hante-Sigalotti} for infinite-dimensional systems, and~\cite{HaidarPepe2020, Haidar-Automatica, HaidarChapter2019, HaidarPepe21} for retarded systems). 

In this paper we give a collection of converse Lyapunov theorems for switching retarded systems with measurable switching signals. The novelty of the obtained results lies in the relaxation of the conditions required by a Lyapunov-Krasovskii functional. We show that the ISS property of a switching retarded system can be characterized by the existence of a continuous (instead of Lipschitz on bounded sets) Lyapunov-Krasovskii functional whose upper right-hand Dini derivative satisfies a dissipation inequality almost everywhere. An important technical tool on which our arguments are based is the recent equivalence property given in~\cite[Theroem 1]{HaidarPepe21} proving that a switching retarded system is ISS (with measurable inputs  and measurable switching signals) if and only if it is ISS for all piecewise-constant inputs and piecewise-constant switching signals. Recall that, when dealing with a retarded system, the map describing the evolution of the state is simply continuous with respect to time~(see, e.g., \cite[Lemma 2.1]{Hale}). Thus a continuous, or even Lipschitz on bounded sets, Lyapunov--Krasovskii functional $V$ evaluated on the solution of a retarded system will be in general continuous and not absolutely continuous with respect to time. By consequence, the nonpositivity of the upper right-hand Dini derivative of $V$ holding almost everywhere, is not sufficient to conclude about the monotonicity of $V$ along the solutions. Thanks to the equivalence property mentioned above, this problem is overcome by restricting the class of inputs and switching signals to the class of piecewise-constant ones. Indeed, in this case, the nonpositivity of the Dini derivative of $V$ along the solutions holds  everywhere instead of almost everywhere permitting to conclude about its monotonicity (see~\cite{Hagood, MironchenkoIto2016}). 

Another contribution of this paper is through the ISS characterization of switching retarded systems using different derivative notions of Driver's and Dini's types. Driver's type derivative (see, e.g., ~\cite{Driver-62, Pepe-Automatica-2007}), by contrast to Dini's one, is an appropriate definition of the derivative of a Lyapunov–Krasovskii functional that does not involve the solution. In~\cite{Pepe-Automatica-2007} it is shown that Driver and Dini derivatives coincide for locally Lipschitz Lyapunov-Krasovskii functionals. Here we extend this result to switching retarded systems. Furthermore, we show that the existence of a Lipschitz on bounded sets Lyapunov-Krasovskii functional whose Driver derivative satisfies a dissipation inequality (which is equivalent, by~\cite[Theorem 2]{HaidarPepe21}, to ISS) is equivalent to the existence of a continuous Lyapunov functional having its Dini derivative satisfying a dissipation inequality almost everywhere. Other  Dini's type derivative definitions, which are used in the literature of retarded systems (see, e.g., ~\cite{8619545}), are also used in the collection of our converse Lyapunov theorems.

The paper is organized as follows. Section~\ref{sec: not} presents the notation, definitions and 
assumptions in use. The statements and proofs of our main results are presented in Section~\ref{main section}. The obtained results 
are discussed in Section~\ref{sec: dis}

\section{Switching retarded systems}\label{sec: not}
In this section we list the notation, definitions, and the main assumptions in use.

\subsection{Notation}
Throughout the paper, we adopt the following notation: $\mathbb{R}$ denotes the set of real numbers, $\mathbb{R}_{+}$ the set of non-negative real numbers, and $\overline{\mathbb{R}}$ the extended real line. By $(\R^n, \|\cdot\|)$ we denote the $n$-dimensional Euclidean space, where $n$ is a positive integer and $\|\cdot\|$ is the Euclidean norm. Given $r>0$, $B(0,r)$ denotes the closed ball of $(\R^n, \|\cdot\|)$ of center $0$ and radius $r$. 
By ${1}_{I}$ we denote the indicator function of a nonempty subset $I$ of $\R$. 

Given $\Delta>0$, $\C:=(\C([-\Delta,0],\mathbb{R}^n),\|\cdot\|_{\infty})$ denotes the Banach space of continuous 
functions mapping $[-\Delta,0]$ into $\mathbb{R}^n$, where $\|\cdot\|_{\infty}$ is the norm of uniform convergence. 
For a function $x:[-\Delta, b)\to\R^n$, with $0< b\leq +\infty$, for $t\in [0,b)$, $x_t:[-\Delta,0]\rightarrow\mathbb{R}^n$
denotes the history function defined by $x_t(\theta)=x(t+\theta)$, $-\Delta\leq\theta\leq0$.
For a positive real $H$ and given $\phi\in \C$, $\C_{H}(\phi)$ denotes the subset $\{\psi\in \C: \|\phi-\psi\|_{\infty}\leq H\}$. We simply denote $\C_{H}(0)$ by $\C_{H}$. 

A measurable function $u:\mathbb{R}_{+}\to\mathbb{R}^{m}$, $m$ positive integer, is said to be essentially bounded if $ess\sup_{t\geq 0}|u(t)|<+\infty$. We use the symbol $\|\cdot\|_{\infty}$ to indicate the essential supremum norm of an essentially bounded function. For given times $0\leq t_1<t_2$, $u_{[t_1,t_2)}:\mathbb{R}_+\to \mathbb{R}^m$ indicates the function given by $u_{[t_1,t_2)}=u(t){1}_{[t_1,t_2)}(t)$ for $t\geq 0$. 
A function $u:\mathbb{R}_{+}\to\mathbb{R}^{m}$ is said to be locally essentially bounded if, for any $t>0$, $u_{[0,t)}$ is essentially bounded.

A function $\alpha:\mathbb{R}_{+}\to\mathbb{R}_{+}$ is said to be of
class $\mathcal{K}$ if it is continuous, strictly increasing and $\alpha(0)=0$; it is said to be
of class $\mathcal{K}_{\infty}$ if it is of class $\mathcal{K}$ and unbounded. A continuous function
$\beta:\mathbb{R}_{+}\times\mathbb{R}_{+}\to\mathbb{R}_{+}$ is said to be of class $\mathcal{KL}$ if $\beta(\cdot,t)$ is of class
$\mathcal{K}$ for each $t\geq 0$ and, for each $s\geq 0$, $\beta(s,\cdot)$ is nonincreasing and converges to zero as $t$ tends to $+\infty$.

With the symbol $\|\cdot\|_{a}$ we indicate any semi-norm in $\C$ such that, for some positive constants $\underline{\gamma_a}$ and $\overline{\gamma_a}$, the following inequalities hold:
\begin{equation*}\label{semi equivalent norm}
\underline{\gamma_a}|\phi(0)|\leq \|\phi\|_{a}\leq \overline{\gamma_a}\|\phi\|_{\infty}, \quad \forall\, \phi\in \C. 
\end{equation*}


\subsection{Definitions and assumptions}
Let us consider the switching control system described by the following retarded functional differential equation 
\begin{equation*}
\S: \begin{array}{llll}
\dot x(t)&=&f_{\sigma(t)}(x_t, u(t)), \quad  & a.e. ~t\geq 0,\\
x(\theta)&=&x_0(\theta), \quad &\theta\in [-\Delta,0], 
\end{array}
\end{equation*}
where: $x(t)\in \R^n$; $n$ is a positive integer; $\Delta$ is a positive real (the maximum involved time delay); 
$x_0\in \C$ is the initial state; the function $\s:\R_{+}\to \mathrm{S}$ is the switching signal; $\mathrm{S}$ is a nonempty set; 
$u:\R_{+}\to\R^m$, $m$ positive integer, is a Lebesgue measurable locally essentially bounded input signal.  

We introduce the following two assumptions:

\begin{assumption}\label{Ass-LBS}
For each $s\in \mathrm{S}$, $f_{s}(0,0)=0$.  Moreover, $f_{s}(\cdot,\cdot)$ is uniformly (with respect to $s\in \mathrm{S}$)
Lipschitz on bounded subsets of $\mathcal{C}\times \R^m$, i.e., for any $H>0$ there exists $L_{H}>0$ such that 
for every $\varphi,\psi\in \C_{H}$ and $u,v\in B(0,H)$, the following inequality holds for all $s\in \mathrm{S}$
\begin{equation*}
|f_s(\varphi,u)-f_s(\psi,v)|\leq L_{H}\left(\|\varphi-\psi\|_{\infty}+|u-v|\right).
\end{equation*}
\end{assumption}

We denote by $\mathcal{U}$ the set of Lebesgue measurable locally essentially bounded inputs from 
$\R_+$ to $\R^m$ and by $\mathcal{U}^{\mathrm{PC}}$ the subset of right-continuous piecewise-constant ones. We denote also by $\mathcal{S}$ the set of measurable signals $\s:\R_+\to \mathrm{S}$ and by $\mathcal{S}^{\mathrm{PC}}$ the subset of right-continuous piecewise-constant ones.

\begin{assumption}\label{Ass-Measurable}
For each $\phi\in \C$, $\s\in \mathcal{S}$ and $u\in {\cal U}$, the function $t\mapsto f_{\s(t)}(\phi,u(t))$, $t\in \mathbb R_+$, is Lebesgue measurable. 
\end{assumption}

Under Assumption~\ref{Ass-LBS} and Assumption~\ref{Ass-Measurable}, the existence and uniqueness of a solution for system~$\S$ as well as its continuous dependence on the initial state is guaranteed by the theory of systems described by retarded functional differential equations (see, e.g., \cite{Hale, Kolmanovskii}). This can be reformulated by the following lemma.

\begin{lemma}\label{existence-uniqueness}
For any  $\phi\in \mathcal{C}$, $u\in\mathcal{U}$ and $\s\in\mathcal{S}$, there exists, uniquely, 
a locally absolutely continuous solution $x(t,\phi,u,\s)$ of $\Sigma$ 
in a maximal time interval $[0,b)$, with $0<b\leq +\infty$. If $b<+\infty$, then the solution is unbounded in $[0,b)$. Moreover, 
for any $\varepsilon>0$, for any $c\in (0,b)$, there exists $\delta>0$ such that, for any $\psi\in \mathcal{C}_{\delta}(\phi)$, the solution 
$x(t,\psi,u,\s)$ exists in $[0,c]$ and, furthermore, the following inequality holds
\begin{equation*}
|x(t,\phi,u,\s)-x(t,\psi,u,\s)|\leq \varepsilon, \quad \forall\,t\in [0,c].
\end{equation*}

\end{lemma}

\medskip
Let us recall the following definition about Driver's form derivative 
of a continuous functional $V:\mathcal{C}\to \mathbb{R}_+$. This definition is a variation 
of the one given in~\cite{Driver-62, Pepe-Automatica-2007,PEPE20061006} for retarded functional differential equations without switching.  
\begin{definition}
For a continuous functional $V:\mathcal{C}\to \mathbb{R}_+$, its Driver's form derivative, 
$D_{(1)}^{+}V:\mathcal{C}\times \R^m\to \overline{\mathbb{R}}$,
is defined, for the switching system $\S$, for $\phi\in \mathcal{C}$ and $u\in \R^m$, as follows,
\begin{equation*}
D_{(1)}^{+}V(\phi,u)=\sup_{s\in \mathrm{S}}\limsup_{h\to0^{+}}\dfrac{V\left(\phi^{\Sigma,s}_{h,u}\right)-V\left(\phi\right)}{h}, 
\end{equation*}
where $\phi^{\Sigma,s}_{h,u}\in \mathcal{C}$ is defined, for $h\in [0,\Delta)$ and $\theta\in [-\Delta,0]$, as follows
\begin{equation*}\label{def-driver}
\phi^{\Sigma,s}_{h,u}(\theta)=\left\{
\begin{array}{lll}
\phi(\theta+h), \quad &\theta\in [-\Delta,-h)\\
\phi(0)+(\theta+h)f_{s}(\phi,u), \quad &\theta\in [-h,0].
\end{array}\right.
\end{equation*}
\end{definition}

Let us also recall the following definition about Dini derivative 
of a continuous functional $V:\mathcal{C}\to \mathbb{R}_+$. This definition is the one given in~\cite{Hale} 
for retarded functional differential equations without switching.

\begin{definition}
Given initial state $\phi\in \mathcal{C}$, $u\in \mathcal{U}$ and $\s\in\mathcal{S}$, for a continuous functional $V:\mathcal{C}\to \mathbb{R}_+$
its Dini derivative $D_{(2)}^{+}V:[0,b)\to \overline{\mathbb{R}}$ is defined, for the switching system $\S$, as follows,
\begin{equation*}
D_{(2)}^{+}V(t)=\limsup_{h\to 0^{+}}\dfrac{V(x_{t+h})-V(x_t)}{h}, 
\end{equation*}
where $x(\cdot)$ is the solution of $\S$ starting from $\phi$ and associated with $u$ and $\s$ over a maximal time interval $[0,b)$. 
\end{definition}

\begin{definition}
For a continuous functional $V:\mathcal{C}\to \mathbb{R}_+$, its $\mathcal{S}$-Dini derivative, 
$D_{(3)}^{+}V:\mathcal{C}\times \mathcal{U}\times\mathcal{S}\to \overline{\mathbb{R}}$,
is defined, for the switching system $\S$, for $\phi\in \mathcal{C}$, $u\in \mathcal{U}$ and $\s\in\mathcal{S}$, as follows,
\begin{equation*}
D_{(3)}^{+}V(\phi,u, \s)=\limsup_{h\to 0^{+}}\dfrac{V(x_{h}(\phi,u,\s))-V(\phi)}{h}, 
\end{equation*}
where $x_h(\phi,u,\s)$ is the solution of $\S$ starting from $\phi$ and associated with $u$ and $\s$. 
\end{definition}

\begin{definition}
For a continuous functional $V:\mathcal{C}\to \mathbb{R}_+$, its mode-Dini derivative, 
$D_{(4)}^{+}V:\mathcal{C}\times \Rset^m\times \mathrm{S}\to \overline{\mathbb{R}}$,
is defined, for the switching system $\S$, for $\phi\in \mathcal{C}$, $v\in \Rset^m$ and $s\in\mathrm{S}$, as follows,
\begin{equation*}
D_{(4)}^{+}V(\phi,v,s)=\limsup_{h\to 0^{+}}\dfrac{V(x_{h}(\phi,v,s))-V(\phi)}{h}, 
\end{equation*}
where $x_h(\phi,v,\s)$ is the solution of $\S$ starting from $\phi$ and associated with $u(t)\equiv v$ and $\s(t)\equiv s$, $t\geq 0$. 
\end{definition}

\begin{definition}
For a continuous functional $V:\mathcal{C}\to \mathbb{R}_+$, its sup-mode-Dini derivative, 
$D_{(5)}^{+}V:\mathcal{C}\times \Rset^m\to \overline{\mathbb{R}}$,
is defined, for the switching system $\S$, for $\phi\in \mathcal{C}$ and $v\in \Rset^m$, as follows,
\begin{equation*}
D_{(5)}^{+}V(\phi,v)=\sup_{s\in\mathrm{S}}D_{(4)}^{+}V(\phi,v,s). 
\end{equation*}
\end{definition}


We give in the following the definition of ISS of system $\S$.

\begin{definition}\label{ISS-def}
We say that system $\S$ is $\mathrm{M}$-$\mathrm{ISS}$ ($\mathrm{PC}$-$\mathrm{ISS}$, respectively), if there exist a function $\beta\in \mathcal{KL}$ and a class $\mathcal{K}$ function $\gamma$ such that, for any $x_0\in \mathcal{C}$, $u\in \mathcal{U}$ ($\mathcal{U}^{\mathrm{PC}}$, respectively) and $\s\in \mathcal{S}$ ($\mathcal{S}^{\mathrm{PC}}$, respectively), the corresponding solution exists in $\mathbb{R}_+$ and, furthermore, satisfies the inequality 
\begin{equation*}
|x(t, x_0,u,\sigma)|\leq \beta(\|x_0\|_{\infty},t)+\gamma(\|u_{[0,t)}\|_{\infty}), \quad\forall\,t\geq 0.
\end{equation*}
\end{definition}

\section{Main results}\label{main section}
The following theorem gives different characterizations of the input-to-state stability property of system $\S$.






\begin{theorem}\label{ISS-converse theorem-Dini} 
The following statements are equivalent: 
\begin{itemize}
    \item [1)] System $\Sigma$ is $\mathrm{PC}$-$\mathrm{ISS}$;
    
     \item [2)] System $\Sigma$ is $\mathrm{M}$-$\mathrm{ISS}$;
     
    \item [3)] there exist a Lipschitz on bounded sets functional $V:\mathcal{C}\to \mathbb{R}_+$, 
functions $\alpha_1, \alpha_2, \alpha_3\in \mathcal{K}_{\infty}$,
and $\alpha_{4}\in \mathcal{K}$ such that the following inequalities hold for every $\phi\in\C$ and $u\in \R^m$:
\begin{enumerate}
\item[(i)]
$\alpha_1(|\phi(0)|)\leq V(\phi)\leq \alpha_2 (\|\phi\|_{a})$,
\item[(ii)]
$D_{(1)}^{+}V(\phi,u)\leq -\alpha_3 (\|\phi\|_{a})+\alpha_{4}(|u|)$;
\end{enumerate}

\item [4)] there exist a continuous functional $V:\mathcal{C}\to \mathbb{R}_+$, 
functions $\alpha_1, \alpha_2, \alpha_3\in \mathcal{K}_{\infty}$,
and $\alpha_{4}\in \mathcal{K}$ such that, for any $\phi\in\C$, any $u\in \mathcal{U}$ and any
$\s\in \mathcal{S}$, the following inequalities hold:
\begin{enumerate}
\item[(i)]
$\alpha_1(|\phi(0)|)\leq V(\phi)\leq \alpha_2 (\|\phi\|_{a}),$
\item[(ii)]
$D_{(3)}^{+}V(x_t,\overline u,\overline\s)\le-\alpha_3(\|x_t\|_a)+\alpha_4(|u(t)|)$, \\$ a.e.\, t\in [0,b)$,\vspace{0.1cm}\\
where $x(\cdot)$ is the solution of $\S$ starting from $\phi$ and associated with $u$ and $\s$ over the maximal interval of definition $[0,b)$, $\overline u(\tau)=u(t+\tau)$ and $\overline\s(\tau)=\s(t+\tau)$, for all $\tau\in [0,b-t)$. \\
Furthermore if $u\in \mathcal{U}^{\mathrm{PC}}$ and $\s\in \mathcal{S}^{\mathrm{PC}}$  then \vspace{0.1cm}
\item[(iii)]
$D_{(3)}^{+}V(x_t,\overline u,\overline\s)\le-\alpha_3(\|x_t\|_a)+\alpha_4(|u(t)|)$, \\ $\forall\, t\in [0,b)$,\vspace{0.1cm}\\
where $x(\cdot)$ is the solution of $\S$ starting from $\phi$ and associated with $u$ and $\s$ over the maximal interval of definition $[0,b)$, $\overline u(\tau)=u(t+\tau)$ and $\overline\s(\tau)=\s(t+\tau)$, for all $\tau\in [0,b-t)$; 
\end{enumerate}

\item [5)] there exist a continuous functional $V:\mathcal{C}\to \mathbb{R}_+$, 
functions $\alpha_1, \alpha_2, \alpha_3\in \mathcal{K}_{\infty}$,
and $\alpha_{4}\in \mathcal{K}$ such that, for any $\phi\in\C$, any $u\in \mathcal{U}$ and any
$\s\in \mathcal{S}$, the following inequalities hold:
\begin{enumerate}
\item[(i)]
$\alpha_1(|\phi(0)|)\leq V(\phi)\leq \alpha_2 (\|\phi\|_{a}),$
\item[(ii)]
$D_{(2)}^{+}V(t)\le-\alpha_3(\|x_t\|_a)+\alpha_4(|u(t)|), \, a.e.\,t\in [0,b)$,\vspace{0.1cm}\\
where $x(\cdot)$ is the solution of $\S$ starting from $\phi$ and associated with $u$ and $\s$ over the maximal interval of definition $[0,b)$.\\
Furthermore if $u\in \mathcal{U}^{\mathrm{PC}}$ and $\s\in \mathcal{S}^{\mathrm{PC}}$  then \vspace{0.1cm}
\item[(iii)]
$D_{(2)}^{+}V(t)\le-\alpha_3(\|x_t\|_a)+\alpha_4(|u(t)|),\, \forall\,t\in [0,b)$, \vspace{0.1cm}
\end{enumerate}
where $x(\cdot)$ is the solution of $\S$ starting from $\phi$ and associated with $u$ and $\s$ over the maximal interval of definition $[0,b)$;

\item [6)] there exist a continuous functional $V:\mathcal{C}\to \mathbb{R}_+$, 
functions $\alpha_1, \alpha_2, \alpha_3\in \mathcal{K}_{\infty}$,
and $\alpha_{4}\in \mathcal{K}$ such that, for any $\phi\in\C$, any $u\in \mathcal{U}^{\mathrm{PC}}$ and any
$\s\in \mathcal{S}^{\mathrm{PC}}$, the following inequalities hold:
\begin{enumerate}
\item[(i)]
$\alpha_1(|\phi(0)|)\leq V(\phi)\leq \alpha_2 (\|\phi\|_{a}),$
\item[(ii)]
$D_{(3)}^{+}V(x_t,\bar u,\bar \s)\le -\alpha_3(\|\phi\|_a)+\alpha_4(|u(t)|),\\ \forall\,t\in [0,b)$, \vspace{0.1cm} 
\end{enumerate}
where $x(\cdot)$ is the solution of $\S$ starting from $\phi$ and associated with $u$ and $\s$ over the maximal interval of definition $[0,b)$, $\overline u(\tau)=u(t+\tau)$ and $\overline\s(\tau)=\s(t+\tau)$, for all $\tau\in [0,b-t)$; 

\item [7)] there exist a continuous functional $V:\mathcal{C}\to \mathbb{R}_+$, 
functions $\alpha_1, \alpha_2, \alpha_3\in \mathcal{K}_{\infty}$,
and $\alpha_{4}\in \mathcal{K}$ such that, for any $\phi\in\C$, any $u\in\Rset^m$ and any
$s\in\mathrm{S}$, the following inequalities hold:
\begin{enumerate}
\item[(i)]
$\alpha_1(|\phi(0)|)\leq V(\phi)\leq \alpha_2 (\|\phi\|_{a}),$
\item[(ii)]
$D_{(4)}^{+}V(\phi,u,s)\le -\alpha_3(\|\phi\|_a)+\alpha_4(|u|)$;
\end{enumerate}

\item [8)] there exist a continuous functional $V:\mathcal{C}\to \mathbb{R}_+$, 
functions $\alpha_1, \alpha_2, \alpha_3\in \mathcal{K}_{\infty}$,
and $\alpha_{4}\in \mathcal{K}$ such that, for any $\phi\in\C$ and any $u\in\Rset^m$ the following inequalities hold:
\begin{enumerate}
\item[(i)]
$\alpha_1(|\phi(0)|)\leq V(\phi)\leq \alpha_2 (\|\phi\|_{a}),$
\item[(ii)]
$D_{(5)}^{+}V(\phi,u)\le -\alpha_3(\|\phi\|_a)+\alpha_4(|u|)$.
\end{enumerate}
\end{itemize}
\end{theorem}

Before giving the proof of Theorem~\ref{ISS-converse theorem-Dini} let us underline what we have mentioned in the introduction concerning the absolute continuity problem of Lyapunov--Krasovskii functionals. In fact, since we deal with a retarded system, the map describing the evolution of the state is simply continuous with respect to time. Thus a continuous (even Lipschitz on bounded sets) Lyapunov--Krasovskii functional evaluated on the solution of such a system will be in general continuous and not absolutely continuous with respect to time (we highlight that this problem is overcome in~\cite{Pepe-Tac-2007} by restricting the class of initial states to continuously differentiable ones; this does not yield any loss of generality because, as it is shown in the same paper, the ISS property holds with continuous initial states if and only if it holds with continuously differentiable ones). By consequence, we cannot directly use the standard comparison lemma~\cite[Lemma 4.4]{LinSontagWang1996} in the proof of the sufficiency parts (i.e., the ones implying the $\mathrm{ISS}$) of Theorem~\ref{ISS-converse theorem-Dini}. Instead, exploiting the equivalence between $\mathrm{M}$-$\mathrm{ISS}$ and $\mathrm{PC}$-$\mathrm{ISS}$ given  by Theorem~\ref{ISS-converse theorem-Dini}, one can use the following comparison lemma from~\cite{MironchenkoIto2016}.

\begin{lemma}{\cite[Lemma 1]{MironchenkoIto2016}}\label{relax_Sontag}
For each continuous and positive definite function $\alpha$, there exists a class $\mathcal{KL}$ function $\beta_{\alpha}$ with the following property: if, for $T>0$ (or $T=+\infty$), $y:[0,T)\to \R_+$ is a continuous non-negative function
which satisfies the inequality 
\begin{equation}\label{beta0}
D^{+}y(t)\leq -\alpha(y(t)), \quad\forall\,t\in[0,T),
\end{equation}
where $D^{+}y$ denotes the upper-right Dini derivative of $y$, with $y(0)=y_0\in \R_+$, then it holds that
\begin{equation}\label{beta1}
y(t)\leq \beta_{\alpha}(y_0,t), \quad\,\forall\, t\in[0,T). 
\end{equation}
\end{lemma}

{\it Proof of Theorem~\ref{ISS-converse theorem-Dini}.}
The proof of $1) \implies 2)$ is given in \cite[Theorem 3.1]{HaidarPepe21}. The proof of $2)\implies 3)$ is given in \cite[Theorem 3.2]{HaidarPepe21}. Concerning the proof of $3)\implies 4)$, let $V$ be the Lipschitz on bounded sets functional given by point 3). Let $x(\cdot)$ be the solution of $\S$ associated with some $\phi\in\C$, $u\in\mathcal{U}$ and $\s\in \mathcal{S}$ over a maximal time interval of definition $[0,b)$. Following the same steps of the proof of~\cite[Theorem 2]{Pepe-Automatica-2007} (see also~\cite{Driver-62}) given for retarded non-switching systems, one can verify that the following equality holds for almost every $t\in [0,b)$
\begin{eqnarray}\label{DD}
   &&\ds\limsup_{h\to 0^{+}}\frac{V(x_{h}(x_t,\bar u,\bar\s))-V(x_t)}{h}\nonumber\\
&&   =\ds\limsup_{h\to0^{+}}\frac{V\left((x_t)^{\Sigma,\s(t)}_{h,u(t)}\right)-V(x_t)}{h}.
\end{eqnarray}
Indeed, Observe that 
\begin{eqnarray*}
&&\ds\frac{V\left((x_t)^{\Sigma,\s(t)}_{h,u(t)}\right)-V\left(x_t\right)}{h}\\
=&&\ds \frac{V\left((x_t)^{\Sigma,\s(t)}_{h,u(t)}\right)-V\left(x_{h}(x_t,\bar u,\bar \s)\right)}{h}\\
&&+\ds\frac{V\left(x_{h}(x_t,\bar u,\bar \s)\right)-V\left(x_t\right)}{h},
\end{eqnarray*}
it is sufficient to prove that for almost every $t\in [0,b)$ we have  
\begin{eqnarray}\label{DD1}
\ds\limsup_{h\to 0^{+}}\dfrac{V\left((x_t)^{\Sigma,\s(t)}_{h,u(t)}\right)-V\left(x_{h}(x_t,\bar u,\bar \s)\right)}{h}=0.
\end{eqnarray}
For this, using the fact that $V$ is Lipschitz on bounded sets, there exists $L=L(x_t)$ such that 
\begin{eqnarray*}
&&\left|V\left((x_t)^{\Sigma,\s(t)}_{h,u(t)}\right)-V\left(x_{h}(x_t,\bar u,\bar \s)\right)\right|\\
&&\hspace{-0.35cm}\leq L\left\|(x_t)^{\Sigma,\s(t)}_{h,u(t)}-x_{h}(x_t,\bar u,\bar \s)\right\|_{\infty}\\
&&\hspace{-0.35cm}=L\sup_{\theta\in [0,h]}\left\|\theta f_{\s(t)}(x_t,u(t))
-\int_{t}^{t+\theta}f_{\s(\tau)}(x_{\tau},u(\tau))d\tau\right\|\\
&&\hspace{-0.35cm}= L\sup_{\theta\in [0,h]}\left\|
\int_{t}^{t+\theta}\left(f_{\s(t)}(x_t,u(t))-f_{\s(\tau)}(x_{\tau},u(\tau))\right)d\tau\right\|\\
&&\hspace{-0.35cm}\leq L\ds\sup_{\theta\in [0,h]}
\int_{t}^{t+\theta}\left\|f_{\s(t)}(x_t,u(t))-f_{\s(\tau)}(x_{\tau},u(\tau))\right\|d\tau\\
&&\hspace{-0.35cm}= L
\int_{t}^{t+h}\left\|f_{\s(t)}(x_t,u(t))-f_{\s(\tau)}(x_{\tau},u(\tau))\right\|d\tau.\\
\end{eqnarray*}
Under Assumption~\ref{Ass-Measurable}, using the Lebesgue's Differentiation Theorem it follows that for almost every $t\in [0,b)$ we have  
\begin{eqnarray*}
\ds\lim_{h\to 0^+}\frac{1}{h}\int_{t}^{t+h}\left|f_{\s(t)}(x_t,u(t)-f_{\s(\tau)}(x_\tau,u(\tau))d\tau)\right|=0.
\end{eqnarray*}
Therefore, equality~\eqref{DD} holds for almost every $t\in [0,b)$. Now observe that 
\begin{eqnarray}\label{DD2}
   &\ds\limsup_{h\to0^{+}}\dfrac{V\left((x_t)^{\Sigma,\s(t)}_{h,u(t)}\right)-V\left(x_t\right)}{h}\leq D_{(1)}^+V(x_t,u(t))\nonumber\\
  & \leq -\alpha_3(\|x_t\|_a)+\alpha_4(|u(t)|).
\end{eqnarray}
From~\eqref{DD} together with~\eqref{DD2} it follows that 
\begin{eqnarray*}
   D^+_{(3)}V(x_t,\bar u,\bar\s)\leq  -\alpha_3(\|x_t\|_a)+\alpha_4(|u(t)|),\, a.e. t\in [0,b).
\end{eqnarray*}
Hence the proof of $3)\implies 4)$. 
Notice that, given any initial state $\phi\in\C$, $u\in\mathcal{U}$ and $\s\in\mathcal{S}$, the following equality holds for all $t\in [0,b)$
\begin{equation}
    D_{(2)}^{+}V(t)=D_{(3)}^{+} V(x_t,\overline u,\overline\s),
\end{equation}
the proof of $4)\implies 5)\implies 6)$ is obvious. The proof of $6) \implies 7)$ follows from the fact that, for each $s\in\mathrm{S}$ and $v\in \Rset^m$, we have $D_{(4)}^{+}V(\phi,v,s)=D_{(3)}^{+}V(\phi,\bar u,\bar\s)$ with $u(\cdot)\equiv v$ and $\s(\cdot)\equiv s$. 
The proof of $7) \implies 8)$ is obvious. Concerning the proof of $8)\implies 1)$, let $\phi\in \mathcal{C}$, $u\in\U^{\mathrm{PC}}$, $\sigma\in \mathcal{S}^{\mathrm{PC}}$, and 
let $x(\cdot)$ be the corresponding solution over a maximal interval of time $[0, b)$, $0<b\leq +\infty$. 
Let $w: [0,b)\to \mathbb{R}_+$ be the function which is defined by
 $$w(t)=V(x_{t}(\phi,u,\sigma)), \quad \forall\,t\in [0,b).$$
Knowing that $u$ and $\sigma$ are piecewise-constants, then for a sufficiently small $h>0$ we have $\s_{|_{[t,t+h)}}\equiv \s(t)$ and $u_{|_{[t,t+h)}}\equiv u(t)$. 
By inequality (ii) of point 8), the following holds for every $t\in [0,b)$
\begin{equation*}
D^{+}w(t)\leq  -\alpha_{3}(\|x_t(\phi,u(t),\s(t))\|_{a})+\alpha_4(|u(t)|).
\end{equation*}
Let the input $u(t)$ be such that $\sup_{t\geq 0}|u(t)|=v$, for a suitable $v\geq 0$. 
By analogous reasoning as in \cite{28018, PEPE20061006}, one can prove the existence of $c\in (0,b]$ such that 
\begin{eqnarray}
&D^{+}w(t)\leq -\alpha(w(t)),\quad &\forall\,t\in [0,c),\label{suf1}\\
&|x(t,\phi,u(t),\s(t))|\leq \gamma(v), \quad 
&\forall\,t\in [c,b),\label{suf2} 
\end{eqnarray}
where $\alpha=\frac{1}{2}\alpha_{3}\circ\alpha_{2}^{-1}$ and $\gamma=\alpha_2\circ\alpha_3^{-1}\circ2\alpha_4$.
Since $t\mapsto w(t)$ is continuous, from Lemma~\ref{relax_Sontag} it holds the existence of a class 
$\mathcal{KL}$ function $\beta_{\alpha}$ such that 
\begin{equation*} 
|w(t)|\leq \beta_{\alpha}(w(0),t), \quad \forall\, t\in [0,c),
\end{equation*} 
from which it follows that 
\begin{equation}\label{ISS1}
|x(t,\phi,u,\s)|\leq \beta(\|\phi\|_{\infty},t), \quad \forall\, t\in [0,c),
\end{equation} 
with $\beta(r,t)=\alpha_{1}^{-1}\circ\beta_{\alpha}(\alpha_2(\overline{\gamma_a}r),t)$. 
By consequence, inequalities~\eqref{suf2} and~\eqref{ISS1} lead to the following inequality 
\begin{equation}\label{ISS}
|x(t,\phi,u,\s)|\leq \beta(\|\phi\|_{\infty},t)+\gamma(v), \quad \forall\, t\in [0,b).
\end{equation} 
It follows, from Lemma~\ref{existence-uniqueness}, that $b=+\infty$.
By causality arguments, and given the arbitrarity of $\varphi\in \C$, $u\in \U^{\mathrm{PC}}$ and $\s\in \mathcal{S}^{\mathrm{PC}}$, the $\mathrm{PC}$-$\mathrm{ISS}$ of system $\S$ is proved. \qed\\

We highlight that, for the cases of nonlinear finite-dimensional and retarded non-switching systems, the result stated in Theorem~\ref{ISS-converse theorem-Dini} concerning the equivalence between items 1) and 2) can be also deduced by~\cite[Theorem 3.3]{Karafyllis2008InputtoOutputSF}, which concerns piecewise-continuous and right-continuous inputs, and by density arguments (see the reasoning used in~\cite[Proposition 3]{Pepe-Tac-2007} for equivalence of ISS with respect to dense sets of initial states).

\section{Conclusions}\label{sec: dis}
In this paper we give a collection of converse Lyapunov theorems for ISS of nonlinear switching retarded systems. In particular, we show that the existence of continuous (instead of locally Lipschitz) Lyapunov-Krasovskii functional whose upper right-hand Dini derivative satisfies a dissipation inequality almost everywhere is necessary and sufficient for the ISS of switching retarded systems. This equivalence property is obtained for a very general class of Lebesgue measurable switching signals. Different derivative notions, which are usually used in the literature of retarded systems, are also used to establish our converse theorems. Future developments may concern  the problem of the input-to-state stabilization and of the input delay tolerance (see, e.g., \cite{wang2020input} and \cite{zhao2022memoryless}) for switching retarded systems.

\bibliography{references}
\end{document}